\documentclass[a4paper,12pt]{article}

\usepackage{amssymb}
\usepackage{amsmath}
\usepackage{hyperref}
\usepackage{amsthm}
\usepackage{amsfonts}
\usepackage[english]{babel}
\usepackage[totalwidth=450pt, totalheight=717pt]{geometry}

\newtheorem{theorem}{Theorem}[section]

\newtheorem{proposition}{Proposition}[section]

\newtheorem{definition}{Definition}[section]

\numberwithin{equation}{section}

\DeclareMathOperator{\dist}   {dist}
\DeclareMathOperator{\dom}    {dom}

\renewcommand{\Im}            {\mathrm{Im}\,}

      \def\dR{{\mathbb R}}

      \def\cF{{\mathcal F}}

      \def\cR{{\mathcal R}}
\def\cS{{\mathcal S}}

\author{Igor Lobanov\thanks{E-mail: lobanov.igor@gmail.com}\\
{\small \emph{ Department of Mathematics}}\\
{\small \emph{St. Petersburg State University of 
		Information Technologies, Mechanics and Optics}}\\
{\small \emph{   Kronverkskiy, 49,  St. Petersburg}}
\and
Vladimir Lotoreichik\thanks{E-mail: vladimir.lotoreichik@gmail.com}\\
{\small \emph{ Department of Mathematics}}\\
{\small \emph{St. Petersburg State University of 
		Information Technologies, Mechanics and Optics}}\\
{\small \emph{   Kronverkskiy, 49,  St. Petersburg}}
 \and
Igor Popov\thanks{E-mail: popov@mail.ifmo.ru}\\
{\small \emph{ Department of Mathematics}}\\
{\small \emph{St. Petersburg State University of 
		Information Technologies, Mechanics and Optics}}\\
{\small \emph{   Kronverkskiy, 49,  St. Petersburg}}
}  

\title{
Lower bound on the spectrum of the Schr\"odinger operator in the plane with $\delta$-potential supported by a curve}

\begin{document}

\maketitle
\begin{abstract}
\noindent
We consider the Schr\"odinger operator in the plane with $\delta$-potential 
supported by a curve $\Gamma$. 
For the cases of an infinite curve and a finite loop 
we give estimates on the lower bound of the spectrum expressed explicitly 
through the strength of the interaction and a parameter $c(\Gamma)$ which characterizes 
geometry of the curve $\Gamma$. Going further we cut the curve into finite number of pieces $\Gamma_i$
and estimate the bottom of the spectrum using parameters $c(\Gamma_i)$.
As an application of the elaborated theory we consider a curve with a 
finite number of cusps and the general ''leaky'' quantum graph. 
\\


\noindent
{\bf Keywords:} 
Schr\"odinger operator, $\delta$-potential, ''leaky'' quantum wires and graphs, 
spectral estimates, cusps, rational curves.\\
\end{abstract}

%
\section{Introduction}
\label{sec:intro}
%
Schr\"odinger operators with $\delta$-potentials on curves 
studied intensively during last two decades are of 
essential interest due to their applications to mesoscopic physics.
Such operators are Hamiltonians of so called ''leaky'' quantum wires 
introduced by P. Exner (for details see the overview \cite{ex:07}).
The operator 
$$H_{\alpha,\Gamma}=-\Delta-\alpha\delta(\cdot-\Gamma),\quad
\dom H\subset L^2(\dR^2).$$
with an attractive $\delta$-potential of a strength $\alpha>0$ having 
a curve $\Gamma$ as the support is strictly defined via the quadratic form
\cite{beks:94}:
$$q_{\alpha,\Gamma}[u]=\|\nabla u\|^2_{L^2(\dR^2)}
-\alpha\|u\|_{L^2(\Gamma)},
\quad \dom q_{\alpha,\Gamma}=H^1(\dR^2).$$
If $\Gamma$ is a straight line,
it is well known that $H_{\alpha,\Gamma}$ has no eigenvalues and
its spectrum is given by  
\begin{equation}
\sigma(H_{\alpha,\Gamma})=\sigma_{ess}(H_{\alpha,\Gamma})= 
\big[-\frac{\alpha^2}{4},+\infty\big),
\end{equation}
assuming that $\alpha$ is a positive constant. 
It is shown in the paper \cite{ei:01}, that if we ''bend'' a little $\Gamma$ 
preserving the smoothness of $\Gamma$, 
then the essential spectrum is preserved and at least one eigenvalue appears. 
For the Schr\"odinger operator with $\delta$-potential on a loop we have according to \cite{beks:94} 
\begin{equation}
\sigma_{ess} = \big[0, +\infty\big),
\end{equation}
and there exists at least one negative eigenvalue.
In the present paper we obtain estimates on the lower bound of 
the spectrum from below for the both cases.

Suppose $\alpha\in\dR_+$, 
$\Gamma$ is a continuous and piecewise-$C^1$ curve in $\dR^2$. 
We need the following characteristic of $\Gamma$:
\begin{equation}\label{eq:c}
c(\Gamma)=\inf_{x,y\in\Gamma}\frac{\dist_{\dR^2}(x,y)}{\dist_\Gamma(x,y)}.
\end{equation}
where we denote by $\dist_{\dR^2}$ distance on the plane and we denote by $\dist_{\Gamma}$  distance along the curve.
It worth noting that if $c(\Gamma)>0$, then $\Gamma$ has no cusps, self-intersections 
and if $\Gamma$ is unbounded then its possible asymptotes are not parallel to each other.
\begin{theorem}
\label{thm:bnd}
Suppose $\Gamma$ is a loop
and $c(\Gamma)>0$.
Then 
\begin{equation}\label{eq:infb}
\inf \sigma(H_{\alpha,\Gamma})\geq -\frac{\alpha^2}{4c(\Gamma)^2}.
\end{equation}
\end{theorem}

\begin{theorem}
\label{thm:unb}
Suppose $\Gamma$ is an infinite curve and $c(\Gamma)>0$.
Then 
\begin{itemize}
\item[(i)] In the most general case we can estimate only the discrete spectrum
\begin{equation}\label{eq:infunb_disc}
\inf \sigma_d(H_{\alpha,\Gamma})\geq -\frac{\alpha^2}{4c(\Gamma)^2}.
\end{equation}
\item[(ii)] But if $\Gamma$ satisfies the condition \eqref{cond:a2} which means the asymptotical straightness then 
\begin{equation}\label{eq:infunb}
\inf \sigma(H_{\alpha,\Gamma})\geq -\frac{\alpha^2}{4c(\Gamma)^2}.
\end{equation}
\end{itemize}
\end{theorem}

The estimate \eqref{eq:infunb} is sharp if $\Gamma$ is a straight line.
If $\Gamma$ consists of two half-lines forming an angle $\beta$,
then $c(\Gamma)=\sin\frac\beta2$, therefore for $\beta\to0$ our 
estimate becomes very rough. However, the estimate \eqref{eq:infunb} can be improved
if $\Gamma$ is a sum of curves with  positive constants $c$ defined in \eqref{eq:c}.

\begin{theorem}
\label{thm:cut}
Suppose the curve $\Gamma$ is a composition of curves $\Gamma_1,\ldots,\Gamma_N$, $N<\infty$, 
in such a way that each $\Gamma_i$ has continuous, piecewise-$C^1$ extension 
$\Gamma_i^{ext} \supset \Gamma_i$ with $c(\Gamma^{ext}_i) >0$.
Moreover we require $\Gamma^{ext}_i$ to be a bounded closed curve or 
an unbounded curve which satisfies \eqref{cond:a2}.
Then the quadratic form $q_{\alpha,\Gamma}$ is closed and semi-bounded from 
below, moreover for the self-adjoint operator $H_{\alpha,\Gamma}$
defined via $q_{\alpha,\Gamma}$ the following estimate holds
\begin{equation}
\label{eq:inf_complex}
\inf\sigma(H_{\alpha,\Gamma})\geq -
N\sum_{i=1}^N\frac{\alpha^2}{4c(\Gamma_i^{ext})^2}.
\end{equation}
\end{theorem}

Consequently, for the mentioned above angle $\beta$ 
we have $\inf\sigma(H_{\alpha,\Gamma})\geq-\alpha^2$.
Moreover, our result can be applied to curves with a finite number of  
cusps, in particular we prove that $H_{\alpha,\Gamma}$ is self-adjoint in this case.

Our strategy for next sections is the following: we give known essential theoretical results, then we prove all the theorems formulated in the introduction. We finish the paper by discussing applications of Theorem \ref{thm:cut} and its possible generalization onto arbitrary ''leaky'' quantum graphs.

\section{Theoretical background}

In this section we recall necessary facts about Schr\"odinger operators with $\delta$-potential on a curve (see \cite{beks:94} and \cite{ei:01} for details).
We consider a Dirac measure $m_{\Gamma}$ associated with the curve $\Gamma$ which is defined as follows,
\begin{equation}
m_\Gamma(M) = l_1(M \cap \Gamma),
\end{equation}
for each Borel set $M \subset \dR^2$, where $l_1$ is the one-dimensional Hausdorff measure; 
in our case it is the length of the part of $\Gamma$ which is inside $M$. 
For $\alpha > 0$ and $c(\Gamma) > 0$ according to~\cite{beks:94} we know that
\begin{equation}
\label{ineq1}
(1 + \alpha) \int_{\dR^2} |\psi(x)|^2 dm_\Gamma(x) \le a \int_{\dR^2} |\nabla \psi(x)|^2dx + b \int_{\dR^2}|\psi(x)|^2dx,
\end{equation}
holds for all $\psi \in \cS(\dR^2)$ and some $a < 1$ and $b$.
We denote by $\cS(\dR^2)$ the Schwartz class of rapidly decreasing functions in $\dR^2$.
The mapping $I_\Gamma \colon \cS(\dR^2) \mapsto L^2(\Gamma)$ defined as $I_\Gamma \psi = \psi$ can be uniquely extended by density argument
\begin{equation}
I_\Gamma \colon W^{1,2}(\dR^2) \mapsto L^2(\Gamma) := L^2(\dR^2, m_\Gamma).
\end{equation}
We use the same symbol to denote $\psi$ as an element of
$L^2(\Gamma)$ and as an element of $L^2(\dR^2)$. 
It is possible to extend the inequality~\eqref{ineq1} 
on $W^{1,2}(\dR^2)$, by exchanging $\psi$ by $I_\Gamma\psi$ in the l.h.s.

The operator we are interested in is defined via the quadratic form 
\begin{equation}
\label{qf}
q_{\alpha, \Gamma}[\psi] := \int_{\dR^2} |\nabla \psi(x)|^2dx - 
\alpha \int_{\dR^2}|I_\Gamma \psi(x)|^2dm_\Gamma(x), 
\end{equation}
with the domain $\dom q_{\alpha, \Gamma} = H^1(\dR^2)$. 
According to the paper \cite{beks:94} under condition~(\ref{ineq1}) this form is closed and bounded from below, with the core $C^{\infty}_0(\dR^2)$. 
This fact implies that a unique self-adjoint operator $H_{\alpha, \Gamma}$ corresponds to this form.

For the operator corresponding to the quadratic form~(\ref{qf}), 
one can apply the Birman-Schwinger technique and obtain the resolvent as follows.  
we denote the resolvent of the Schr\"odinger operator in $L^2(\dR^2)$ without potential 
by $R^k_0$, $\Im k>0$, that is $(-\Delta-k^2)^{-1}=R^k_0$.
The resolvent $R^k_0$ appears to be the integral operator with the kernel
\begin{equation}
\label{green}
G_{k}(x - y) = \frac{i}{4}H_{0}^{(1)}(k|x - y|)
\end{equation}

Let $\mu$ and $\nu$ be arbitrary positive Radon measures in $\dR^2$, satisfying the condition $\mu(x) = \nu(x) = 0$ for each $x \in \dR^2$.
By  $R^k_{\nu,\mu}$ we denote an integral operator from $L^2(\mu) := L^2(\dR^2, \mu)$ to $L^2(\nu)$ with the kernel $G_k$, i.e.
\begin{equation}\label{eq:Rknm}
R^k_{\nu,\mu}\phi = G_k * \phi\mu
\end{equation}
holds a.e. with respect to the measure $\nu$ for all $\phi \in \dom R^k_{\nu,\mu} \subset L^2(\mu)$.
In our case there are two measures: $m_\Gamma$ introduced before and Lebesgue measure $dx$ in $\dR^2$. 
In this notation according to \cite{beks:94} we have the following
\begin{proposition}
\label{krein}
\begin{itemize}
\item[(i)] There exists $\kappa_0 > 0$ such that $I - \alpha R^{i\kappa}_{m_\Gamma,m_\Gamma}$ on  $L^2(\Gamma)$ has the bounded inverse for all $\kappa > \kappa_0$.
\item[(ii)] Let $\Im k > 0$. Assume that  $I - \alpha R^{k}_{m_\Gamma,m_\Gamma}$ is invertible and the operator
\begin{equation}
\label{kreinf}
R^k := R^k_0 + \alpha R^k_{dx,m_\Gamma}\big[I - \alpha R^k_{m_\Gamma,m_\Gamma}\big]^{-1}R^k_{m,dx}
\end{equation}
from $L^2(\dR^2)$ to $L^2(\dR^2)$ is everywhere defined. Then $k^2$ belongs to $\rho(H_{\alpha, \Gamma})$ and $\big(H_{\alpha, \Gamma} - k^2\big)^{-1} = R^k$.
\item[(iii)]
$\dim \ker\big(H_{\alpha, \Gamma} - k^2\big) = \dim \ker\big(I -\alpha R^k_{m_\Gamma,m_\Gamma}\big)$ for all $k$ such that $\Im k > 0$.
\end{itemize}
\end{proposition}

Next we give some known results about the the spectrum of $H_{\alpha, \Gamma}$ in cases of a closed loop and of an unbounded curve.
We assume that the curve $\Gamma$ is defined by continuous piecewise-$C^1$ mapping 
$\gamma = (\gamma_x,\gamma_y)$, which maps from $\dR$ to $\dR^2$ in the case of an unbounded curve and maps
from $[0,L]$ to $\dR^2$ with $\gamma(L) = \gamma(0)$ in the case of a loop. 
We assume that the curve is defined in a natural parametrization, i.e.  $|\dot \gamma(s)| = 1$. 

Since $\gamma$ is continuous and is naturally parameterized, we have obviously the following property
\begin{equation}
\label{ineq2}
|\gamma(s) - \gamma(s')| \le |s - s'|,
\end{equation}
for all $s, s' \in \dR$. 
Next we give assumptions which were introduced in \cite{ei:01}
\begin{itemize}
\item[(a1)]\label{cond:a1} The constant $c(\Gamma)$ is strictly positive .
In particular for an infinite curve it means that 
$$|\gamma(s) - \gamma(s')| \ge c(\Gamma)|s - s'|,$$
and for a finite loop it means that
$$|\gamma(s) - \gamma(s')| \ge c(\Gamma) \min\{|s-s'|, L - |s-s'|\}.$$
\item[(a2)] This assumption is meaningful only for an unbounded curve. We require from $\Gamma$ to be asymptotically straight in the following sense: there exist positive constants $d, \mu > \frac12$ and $\omega \in (0, 1)$ such that
\begin{equation}
\label{cond:a2}
1 - \frac{|\gamma(s) -\gamma(s')|}{|s - s'|} \le d\big[1 + |s + s'|^{2\mu}\big]^{-\frac12}
\end{equation}
holds in the sector $S_{\omega} =\big\{(s,s') \colon \omega \le \frac{s}{s'} \le \omega^{-1}\big\}$.
\end{itemize}

Now we can formulate some of the results from papers \cite{beks:94} and \cite{ei:01} which were non-formally given in the introduction.
\begin{theorem}
\label{thm:exner}
 Assume $\Gamma$ is continuous, piecewise-$C^1$ curve and (a1) holds then
\begin{itemize}
 \item [(i)] If $\Gamma$ is a loop then 
\begin{equation}
\sigma_{ess}(H_{\alpha, \Gamma}) = \big[0, +\infty\big)
 \end{equation}
and at least one negative eigenvalue exists.
\item[(ii)] If $\Gamma$ is an unbounded curve and in addition (a2) holds then
\begin{equation}
\sigma_{ess}(H_{\alpha, \Gamma}) = \big[ -\frac{\alpha^2}{4},+\infty\big)
\end{equation}
Moreover if $\Gamma$ is not a straight line then at least one eigenvalue below the essential spectrum exists.
\end{itemize}
\end{theorem}

\section{Proof of  Theorem \ref{thm:bnd}}
\label{sec:proof:bnd}
\begin{proof}
We introduce on $L^2[0,L]$ an integral operator  
$\cR^{\kappa}_{\alpha, \Gamma} := \alpha R^{i\kappa}_{m_\Gamma ,m_\Gamma}$,
where $R^{\kappa}_{\nu\,\mu}$ is defined by \eqref{eq:Rknm}. 
According to Proposition~\ref{krein} item (iii), 
a negative value $\lambda = - \kappa^2$ is an eigenvalue of the operator 
$H_{\alpha, \Gamma}$ if and only if  the equation 
\begin{equation}
\label{int_eq_bnd}
\cR^{\kappa}_{\alpha, \Gamma} \psi = \psi,
\end{equation}
has non-trivial solution.
The integral kernel of $\cR^{\kappa}_{\alpha, \Gamma}$ has the following form
\begin{equation}
\cR^{\kappa}_{\alpha, \Gamma}(s,s') = \frac{\alpha}{2\pi}K_0(\kappa|\gamma(s) - \gamma(s')|),
\quad s,s' \in [0, L],
\end{equation} 
where $K_0$ is the modified Bessel function of the 2$^{nd}$ kind (Macdonald function); recall $K_0(z) = \frac{\pi i}{2}H_0^{(1)}(iz)$ \cite{as:64}.
The equation (\ref{int_eq_bnd}) has only trivial solution if  $\|\cR^{\kappa}_{\alpha, \Gamma}\| < 1$. 

From the property (a1) and monotonous decreasing of the Macdonald function it follows that
\begin{equation}
\label{ineq3_bnd}
0 \le K_0(\kappa|\gamma(s) - \gamma(s')|) \le K_0(\kappa c(\Gamma)\dist_{\Gamma}(\gamma(s),\gamma(s'))).
\end{equation}
We introduce the operator $\cR^{\kappa}_c$ in $L^2(\Gamma)$ with the integral kernel 
\begin{equation*}
\cR^{\kappa}_c(s,s') = \frac{\alpha}{2\pi}K_0(\kappa c(\Gamma)\dist_{\Gamma}(\gamma(s), \gamma(s')))=\frac{\alpha}{2\pi}K_0(\kappa c(\Gamma)\min\{|s-s'|, L - |s-s'|\}).
\end{equation*}
According to~(\ref{ineq3_bnd})
\begin{equation}
\|\cR^{\kappa}_{\alpha, \Gamma}\| \le \|\cR^{\kappa}_c\|.
\end{equation}
Since the integral kernel of $\cR^{\kappa}_c$ is positive, we can
apply Schur's Lemma (Theorem 5.2 in \cite{hs:78} ) to estimate its norm from above
\begin{equation}
\label{Schur}
\|\cR^{\kappa}_c\| \le \sqrt{\sup_{s \in [0,L]} \int_0^{L}\cR^{\kappa}_c(s,s')ds'\cdot\sup_{s' \in [0,L]} \int_0^{L} \cR^{\kappa}_c(s,s')ds} 
\end{equation}
Since the operator $\cR^{\kappa}_c$ has symmetric kernel then both supremums are equal.
Assume $s < \frac{L}{2}$, then
\begin{eqnarray*}
\int_0^{L}\cR^{\kappa}_c(s,s')ds' &=&\frac{\alpha}{2\pi}\int_0^{L}
K_0(\kappa c(\Gamma)\min\{|s-s'|, L - |s-s'|\})ds'=\\
&=& \frac{\alpha}{2\pi} \int_{s-L}^s K_0(\kappa c(\Gamma)\min\{|t|, L - |t|\})dt =\\
&=&\frac{\alpha}{2\pi} \Big(\int_{s-L}^{-\frac{L}{2}} K_0(\kappa c(\Gamma)(L - |t|))dt + 
\int_{-\frac{L}{2}}^{s} K_0(\kappa c(\Gamma)|t|))dt\Big) =\\
&=&\frac{\alpha}{2\pi} \Big(\int_{s}^{\frac{L}{2}} K_0(\kappa c(\Gamma)t)dt + 
\int_{-\frac{L}{2}}^{s} K_0(\kappa c(\Gamma)|t|))dt\Big) = \\
&=& \frac{\alpha}{2\pi} \int_{-\frac{L}{2}}^{\frac{L}{2}} K_0(\kappa c(\Gamma) |t|)dt < \frac{\alpha}{2\pi} \int_{-\infty}^{+\infty} K_0(\kappa c(\Gamma) |t|)dt = \frac{\alpha}{2c(\Gamma)\kappa}
\end{eqnarray*}
The same estimate holds for $s > \frac{L}{2}$ also. 
Hence, $\|\cR^{\kappa}_{\alpha, \Gamma}\| < \frac{\alpha}{2c(\Gamma)\kappa}$ and if $\kappa >\frac{\alpha}{2c(\Gamma)}$
then the r.h.s. of the last estimate is less than one. 
Therefore if $\lambda < -\frac{\alpha^2}{4c(\Gamma)^2}$ then $\lambda$ can not be an eigenvalue. Combining this statement with the item (i) of Theorem~\ref{thm:exner} we obtain
\begin{equation*}
\inf \sigma(H_{\alpha, \Gamma}) \ge -\frac{\alpha^2}{4c(\Gamma)^2}
\end{equation*}

\end{proof}

\section{Proof of  Theorem \ref{thm:unb}}
\begin{proof}
The first part of this proof is the repeating of the previous one with small changes.
We introduce on $L^2(\dR)$ an integral operator  
$\cR^{\kappa}_{\alpha, \Gamma} := \alpha R^{i\kappa}_{m_\Gamma ,m_\Gamma}$
where $R^{\kappa}_{\nu,\mu}$ is defined by \eqref{eq:Rknm}. 
According to Proposition~\ref{krein} item (iii), 
a negative value $\lambda = - \kappa^2$ is an eigenvalue of the operator 
$H_{\alpha, \Gamma}$ if and only if  the equation 
\begin{equation}
\label{int_eq}
\cR^{\kappa}_{\alpha, \Gamma} \psi = \psi,
\end{equation}
has non-trivial solution.
The integral kernel of $\cR^{\kappa}_{\alpha, \Gamma}$ has the following form
\begin{equation}
\cR^{\kappa}_{\alpha, \Gamma}(s,s') = \frac{\alpha}{2\pi}K_0(\kappa|\gamma(s) - \gamma(s')|),
\end{equation} 
where $s, s' \in \dR$.
The equation (\ref{int_eq}) has only trivial solution if  $\|\cR^{\kappa}_{\alpha, \Gamma}\| < 1$. 
From the property (a1) and monotonous decreasing of the Macdonald function it follows that
\begin{equation}
\label{ineq3}
0 \le K_0(\kappa|\gamma(s) - \gamma(s')|) \le K_0(\kappa c(\Gamma)|s- s'|).
\end{equation}
We introduce an operator $\cR^{\kappa}_c$ with the integral kernel $\cR^{\kappa}_c(s,s') = \frac{\alpha}{2\pi}K_0(\kappa c(\Gamma)|s- s'|)$.
According to~(\ref{ineq3})
\begin{equation}
\|\cR^{\kappa}_{\alpha, \Gamma}\| \le \|\cR^{\kappa}_c\|.
\end{equation}
In the case of an unbounded curve we can be more precise than in the bounded case. 
We calculate the norm of $\cR^{\kappa}_c$ exactly. 
Let $\cF \colon L^2(\dR) \to L^2(\dR)$ be the Fourier transform on the real axis defined as in \cite{be:69} by the following formula
\begin{equation*}
\hat f(p) = (\cF f)(p) = \frac{1}{\sqrt{2\pi}}\int_{\dR} f(\xi)e^{-i \xi p}d\xi
\end{equation*}
with the inverse transform
\begin{equation*}
\check f(p) = (\cF^{-1} f)(p) = \frac{1}{\sqrt{2\pi}}\int_{\dR} f(\xi)e^{i \xi p}d\xi.
\end{equation*}
Next we consider the operator  $\hat \cR^{\kappa}_c := \cF \cR^{\kappa}_c \cF^{-1}$ (momentum representation).  
Its norm coincides with the norm of $\cR^{\kappa}_c$, since the Fourier transform is unitary, furthermore
\begin{equation*}
\cF \cR^{\kappa}_c \cF^{-1} f = \frac{\alpha}{2\pi}\cF\big( K_0(\kappa c|\xi|) * (\cF^{-1} f)\big) = \frac{\alpha}{\sqrt{2\pi}} \cF (K_0(\kappa c(\Gamma) |\xi|)) \cF \cF^{-1}f =  \frac{\alpha}{2\sqrt{c(\Gamma)^2\kappa^2 + p^2}} f,
\end{equation*}
here we used the formula $\cF (f * g) = \sqrt{2\pi}\hat f \hat g$ and the Fourier transform of the Macdonald function (see \cite{be:69}).
The operator $\hat \cR^{\kappa}_c$ is the multiplication operator on the function $\frac{\alpha}{2\sqrt{c(\Gamma)^2\kappa^2 + p^2}}$ of the argument $p$ in the space $L^2(\dR)$, hence
\begin{equation*}
\|\hat \cR^{\kappa}_c\| = \sup_{p \in \dR}~ \frac{\alpha}{2\sqrt{c(\Gamma)^2\kappa^2 + p^2}} = \frac{\alpha}{2c(\Gamma)\kappa}.
\end{equation*}
This norm is less than one when $\kappa > \frac{\alpha}{2c(\Gamma)}$. 
Hence, if $\lambda < -\frac{\alpha^2}{4c(\Gamma)^2}$ then it can not be an eigenvalue of $H_{\alpha, \Gamma}$ and the first item is proven.
In order to prove item (ii) we use the second statement of Theorem \ref{thm:exner} and we obtain
\begin{equation*}
 \inf \sigma(H_{\alpha, \Gamma}) \ge -\frac{\alpha^2}{4c(\Gamma)^2}
\end{equation*}

\end{proof}

\section{Proof of Theorem \ref{thm:cut}}
\begin{proof}
First of all we prove self-adjointness of $H_{\alpha, \Gamma}$.
For each $\Gamma_i$ we know that $c(\Gamma_i) > 0$ then
we can write down the inequality of the type~\eqref{ineq1},
\begin{equation}
N(1 + \alpha)\|u\|_{L^2(\Gamma_i)}^2 \le a_i \|\nabla u\|_{L^2(\dR^2)}^2 + b_i\|u\|_{L^2(\dR^2)}^2,
\end{equation}
where $0 < a_i <1$ and $b > 0$,
since \eqref{ineq1} is valid for every $\alpha\in\dR$.
Adding these inequalities for all $\Gamma_i$ and then dividing by $N$ we obtain
\begin{equation*}
 (1 + \alpha)\|u\|_{L^2(\Gamma)} \le \frac{\sum_{i=1}^N a_i}{N} \|\nabla u\|_{L^2(\dR^2)}^2 +  \frac{\sum_{i=1}^N b_i}{N}\|u\|_{L^2(\dR^2)}^2.
\end{equation*}
Since $0 <\frac{\sum_{i=1}^N a_i}{N} < 1 $ then according to \cite{beks:94} the form $q_{\alpha, \Gamma}$ with $H^1(\dR^2)$ as domain is bounded from below and closed.
A unique self-adjoint operator $H_{\alpha, \Gamma}$ corresponds to this form.
Next we prove the estimate~\eqref{eq:inf_complex}.
We consider a quadratic form for each $\Gamma_i$
\begin{equation}
\label{qf_i}
q_{\alpha N, \Gamma_i}[u] = \|\nabla u\|^2 - \alpha N\|u\|_{L^2(\Gamma_i)}^2,
\end{equation}
with $\dom q_{\alpha N, \Gamma_i} = H^1(\dR^2)$. 
We consider also quadratic forms for each $\Gamma_i^{ext}$  
\begin{equation*}
q_{\alpha N, \Gamma_i^{ext}}[u] = \|\nabla u\|^2 - \alpha N\|u\|_{L^2(\Gamma_i^{ext})}^2,
\end{equation*}
Since $\Gamma_i^{ext} \supset \Gamma_i$ then
\begin{equation*}
q_{\alpha N, \Gamma_i}[u] \ge q_{\alpha N, \Gamma_i^{ext}}[u],\quad \forall u \in H^1(\dR^2).
\end{equation*}
The quadratic form for the Schr\"odinger operator with $\delta$-potential of strength $\alpha$ on the whole $\Gamma$ can be expressed as follows
\begin{equation}
q_{\alpha, \Gamma}[u] = \frac1N\sum_{i=1}^N q_{\alpha N, \Gamma_i}[u] \ge \frac1N\sum_{i=1}^N q_{\alpha N, \Gamma_i^{ext}}[u] 
\end{equation}
The lower bound of the form $q_{\alpha N, \Gamma_i^{ext}}$ coincides with $\inf \sigma(H_{\alpha N, \Gamma_i^{ext}})$.
If $\Gamma_i^{ext}$ is a loop we apply Theorem~\ref{thm:bnd}. If it is unbounded we apply item (ii) of Theorem~\ref{thm:unb}. 
Since, $\inf \sigma(H_{\alpha, \Gamma}) = \inf_{u \in H^1(\dR^2)} q_{\alpha, \Gamma}[u]$ then
\begin{equation}
\inf \sigma(H_{\alpha, \Gamma}) \ge -
N\sum_{i=1}^N\frac{\alpha^2}{4c(\Gamma_i^{ext})^2}
\end{equation}
\end{proof}

\section{Applications and possible generalizations}

\subsection{Application to curves with cusps}

First of all we give the definition of a cusp (see \cite{po:94}).
\begin{definition}
\label{cusp}
For the curve defined as a zero set of a function of two variables $F(x,y) \in C^{\infty}$ the point $(x_0, y_0)$ is a cusp if 
\begin{itemize}
 \item[(i)] $F(x_0,y_0) = 0$.
\item[(ii)] $\frac{\partial F}{\partial x}(x_0,y_0)  =\frac{\partial F}{\partial y}(x_0,y_0) = 0$.
\item[(iii)]  The Hessian matrix of second derivatives has zero determinant at the point $(x_0,y_0)$. That is:
$$\begin{vmatrix} \frac{\partial^2 F}{\partial x^2} & \frac{\partial^2 F}{\partial x\,\partial y} \\ \\ \frac{\partial^2 F}{\partial
 x\,\partial y} & \frac{\partial^2 F}{\partial y^2} \end{vmatrix}(x_0,y_0) =0.$$
\end{itemize}
\end{definition}

Next we give examples of rational curves with cusps 
\begin{itemize}
\item[(i)] $\gamma(t) = (t^2,t^3)$ is a spinode.
\item[(ii)]$\gamma(t) = (t^2, t^5)$ is a rhamphoid.
\item[(iii)]$\gamma(t)=\gamma_{m, n}(t) := (t^{2n}, t^{2m+1})$, where $m \ge n$.
\end{itemize}
All these curves have a cusp at the origin which can be simply checked by constructing appropriate $F(x,y)$ and verifying Definition \ref{cusp}.
Let us show for spinode that $c(\Gamma) = 0$.
\begin{eqnarray*}
 c(\Gamma) &\le& \lim_{t_0 \to 0+}\frac{\sqrt{(t_0^2 - t_0^2)^2+ (t_0^3 + t_0^3)^2}}{\int_{-t_0}^{t_0} \sqrt{4t^2 + 9t^4}dt } = 
= \lim_{t_0 \to 0+} \frac{2 t_0^3}{\frac{16}{27}\big((1 + \frac{9}{4}t_0^2\big))^\frac32 - 1\big)}\\
& =& \lim_{t_0 \to 0+} \frac{2t_0^3}{2t_0^2 +o(t_0^2)} = 0
\end{eqnarray*}

Consider the curve $\Gamma$ without self-intersections. We assume that $\Gamma$ is defined by mapping $\gamma(t)$
satisfying the condition (a2). We also assume that $\Gamma$ has finitely many cusps  at points
$t_1 < t_2 <... < t_N$. Next, we accomplish this sequence with points $t_0 = -\infty$ and $t_{N+1} = +\infty$.
We make a cut of the curve $\Gamma$ into parts $\Gamma_i = \gamma([t_i, t_{i+1}]), ~i = 0..N$. If each part has an extension $\Gamma_i^{ext}$ such that $c(\Gamma_i^{ext}) > 0$, then we can apply Theorem \ref{thm:cut} to prove self-adjointness of $H_{\alpha, \Gamma}$ and to obtain the estimate on the spectrum from below.

\subsection{Generalization to arbitrary ''leaky'' quantum graphs}

The result of Theorem \ref{thm:cut} can be extended onto arbitrary ''leaky'' quantum graphs.
Recall that the Hamiltonian of ''leaky'' quantum graphs is defined by 
the self-adjoint operator $H_{\alpha,\Gamma}$ corresponding to the quadratic form  $q_{\alpha,\Gamma}$
given by \eqref{qf} where $\Gamma$ is a metric graph (a union of segments
intersecting only at ends).
We assume that the graph is finite (i.e. it consists of finite number of segments)
and all infinite segments are asymptotically straight.
Then the Theorem \ref{thm:cut} can be reformulated
\begin{theorem}
Assume that for every segment $\Gamma_i$ of a metric graph $\Gamma$
we can construct an extension $\Gamma_i^{ext} \supset \Gamma_i$ which is an infinite curve satisfying \eqref{cond:a2} or a loop, continuous and piecewise-$C^1$ with 
$c(\Gamma_i^{ext})>0$. 
Then the quadratic form $q_{\alpha,\Gamma}$ is closed and bounded from below and a unique self-adjoint operator $H_{\alpha, \Gamma}$  corresponds to it. 
Moreover the following estimate is valid
\begin{equation}
 \inf \sigma(H_{\alpha, \Gamma}) \ge -N\sum_{i=1}^{N} \frac{\alpha^2}{4c(\Gamma_i^{ext})^2},
\end{equation}
where $N$ is the number of segments.
\end{theorem}

\section{Acknowledgements}
This work was partially supported by the grant RFBR-NANU no. 09-01-90410.
\label{sec:est}

\newcommand{\etalchar}[1]{$^{#1}$}
\providecommand{\bysame}{\leavevmode\hbox to3em{\hrulefill}\thinspace}
\providecommand{\MR}{\relax\ifhmode\unskip\space\fi MR }
\renewcommand{\MR}[1]{}  


\begin{thebibliography}{EKK{\etalchar{+}}08}


\bibitem[AS64]{as:64}
M.~Abramovitz, I.~Stegun, \emph{Handbook of mathematical functions with formulas, graphs and mathematical tables}, National Bureau of Standards, Applied Mathematics Series, 55, 1972.

\bibitem[BE69]{be:69}
H. Bateman, A. Erdelyi,
\emph{Higher Transcendental Functions}, McGraw-Hill, 1953.


\bibitem[BEKS94]{beks:94}
J.~F.~Brasche, P.~Exner, Yu.~A.~Kuperin, P.~Seba,
\emph{Schr\"odinger operators with singular interactions}
J. Math. Anal. Appl. \textbf{184} (1994), 112-139.


\bibitem[BEW08]{bew:08}
B.~M.~Brown, M.S.P.~Eastham, I.~G.~Wood,
\emph{Conditions for the spectrum associated with an asymptotically straight leaky wire to contain an interval} $[- \frac{\alpha^2}{4}, +\infty)$, Archiv der Math. \textbf{90}, no. 6,	554-558.

\bibitem[Ex06]{ex:06}
P.~Exner,
\emph{An isoperimetric problem for leaky loops and related mean-chord inequalities} J. Math. Phys. {\bf 46} (2005), 062105. 

\bibitem[Ex07]{ex:07}
P.~Exner,
\emph{Leaky quantum graphs: a review },
in Proceedings of the INI programme "Analysis on Graphs and Applications" (Cambridge 2007), AMS "Proceedings of Symposia in Pure Mathematics" Series, vol. \textbf{77}, Providence, R.I., (2008); pp. 523-564.

\bibitem[EN01]{en:01}
P.~Exner, K.~N\v{e}mcov\'{a}: {\em Bound states in point
interaction star graphs}, J. Phys. A: Math. Gen. {\bf 34} (2001), 7783--7794.

\bibitem[EN03]{en:03}
P.~Exner, K.~N\v{e}mcov\'{a}: {\em Leaky quantum graphs: approximations by point interaction Hamiltonians
}, 
J. Phys. A {\bf 36} (2003), 10173-10193

\bibitem[EI01]{ei:01}
P.~Exner and T.~Ichinose,
\emph{Geometrically induced spectrum in curved leaky wires},
J. Phys. A\textbf{34},
(2001), 1439--1450.

\bibitem[EK03]{ek:03}
P.~Exner, S.~Kondej,
\emph{Bound states due to a strong delta interaction supported by a curved surface},
J. Phys. A\textbf{36} (2003), 443-457

\bibitem [EY02]{ey:02} P.~Exner, K.~Yoshitomi, \emph{Asymptotics of
eigenvalues of the Schr\"odinger operator with a strong
$\delta$-interaction on a loop}, J. Geom. Phys. \textbf{41}
(2002), 344--358.

\bibitem[KMR97]{kmr:97} V.~Kozlov, , V.~Maz'ya, J.~Rossmann  
\emph{Elliptic boundary value problems in domains with point singularities}, Mathematical Surveys and Monographs {\bf 52} (1997), Amer. Math. Soc.

\bibitem[MNP00]{mnp:00} V.~G.~Maz'ya,  Yu.~V.~Netrusov, S.~V.~Poborchiy,
\emph{Boundary values of functions in Sobolev spaces on certain non-Lipschitzian domains}
St. Petersburg Mathematical Journal, 2000, {\bf 11}:1, 107-128. 

\bibitem[HS78]{hs:78} P. R. Halmos,  V.~S.~Sunder, \emph{Bounded integral operators on L2 spaces}, 
Ergebnisse der Mathematik und ihrer Grenzgebiete (Results in Mathematics and Related Areas), vol. 96. Springer-Verlag, Berlin, 1978.

\bibitem[Po94]{po:94} 
I.~Porteous,
\emph{Geometric Differentiation - for the intelligence of curves and surfaces},
Cambridge University Press", 1994.



\bibitem[Te90]{teta:90}
A.~Teta,
\emph{Quadratic forms for singular perturbations of the Laplacian},
Res. Inst. Math. Sci. \textbf{26},
(1990), no. 5, 803--817.


\end{thebibliography}
\end{document}